\documentclass[a4paper,11pt]{article}
\usepackage[utf8]{inputenc}
\usepackage[T1]{fontenc}
\usepackage[english]{babel}

\usepackage[tbtags]{amsmath}
\usepackage{amssymb}
\usepackage{amsfonts}
\usepackage{amsthm}
\usepackage{calrsfs}
\usepackage{array}
 \usepackage{bbm}
\usepackage{stmaryrd}
\usepackage{upgreek}
\usepackage[svgnames]{xcolor}
\usepackage{hyperref}

\def \p {\psi (y,s)}
\numberwithin{equation}{section}
\def\iint {\int_{B}}
\def\b{{\frac{1}{s^b}}}
\def\grad{\nabla}
\def\div{\, \mbox{div}\,  }

\def\ia{\int_{s}^{s+1}}

\def\L{L(w(s),s)}
\def\q{\phi (y,s)}

\def\Box{\hfill\rule{2.5mm}{2.5mm}}
\def\t{{\mathrm{d}}\tau}

\def\y{{\mathrm{d}}y}
\def\e{\varepsilon}
\def\no{\nonumber}

\def \p {\psi (y,s)}
\def\grad{\nabla}
\def\div{\, \mbox{div}\,  }

\def \er {\mathbb R}

\def\ta{\phi (y,\tau)}
\def\R{{\mathbb {R}}}

\def\cprime{$'$}

\newcommand{\ds}{\displaystyle}

\def\t{{\mathrm{d}}\tau}

\def\y{{\mathrm{d}}y}
\def\e{\varepsilon}
\def\no{\nonumber}

\def\y{{\mathrm{d}}y}

\theoremstyle{plain}
\newtheorem{thm}{Theorem}[section]
\newtheorem*{thm*}{Theorem}
\newtheorem{pro}[thm]{Proposition}

\newtheorem{cor}[thm]{Corollary}
\newtheorem{lem}[thm]{Lemma}

\theoremstyle{definition}

\theoremstyle{remark}

\makeatletter
\def\blfootnote{\xdef\@thefnmark{}\@footnotetext}
\makeatother

\title{\bf The blow-up rate for strongly perturbed semilinear wave equations  in the conformal regime
without a radial assumption}

\author {M.A. Hamza\footnote{The author is partially supported by the ERC Advanced Grant no.291214, BLOWDISOL during his visit to LAGA, Universit\'{e} Paris 13 in 2014.  }\\
{\it \small Facult\'e des Sciences de Tunis, Universit\'e Tunis El Manar, Tunis, Tunisie}}

\allowdisplaybreaks

\begin{document} 

\date{}
\maketitle
\begin{abstract}
We consider in this paper a large class of perturbed semilinear wave
equations 
with critical (in the conformal transform sense) power nonlinearity.  We will show 
 that the blow-up rate of any singular solution is given by the solution of the non-perturbed associated ODE.
The result in the radial case  has been proved in  \cite{omar2}.  
The same approach will be followed
here, but the main difference 
 is to construct a Lyapunov functional
in similarity variables valid
in  the non-radial case,
 which is far from being trivial. That functional is obtained 
by combining some classical     estimates 
and a new identity of the Pohozaev type  
 obtained  by multiplying the equation \eqref{A}  by $y.\grad w$ in  a suitable weighted space. 
\end{abstract}

\medskip



\section{Introduction}
This paper is devoted to the study of blow-up solutions for the
following semilinear  wave equation:
\begin{equation}\label{gen}
\left\{
\begin{array}{l}
\partial_t^2 u =\Delta u+|u|^{p-1}u+f(u)+g(x,t,\nabla u,\partial_t u ),\\
\\
u(x,0)=u_0(x)\in  H^{1}_{loc,u}(\er^N)\qquad {\textrm{and}}\quad
\partial_tu(x,0)=u_1(x)\in  L^{2}_{loc,u}(\er^N),\\
\end{array}
\right.
\end{equation}
with conformal power nonlinearity
\begin{equation}\label{pc}
p=p_c\equiv 1+\frac4{N-1}, \quad {\textrm{where}} \quad N\ge 2,
\end{equation}
and $u(t):x\in \er^N \rightarrow u(x,t)\in \er$. The space $L^2_{loc,u}$
is the set of all $v\in L^2_{loc}$ such that
$$\|v\|_{L^2_{loc,u}}:= \sup_{d\in{\er^N}}\Big(\int_{|x-d|<1} |v(x)|^2
 {\mathrm{d}}x \Big)^{\frac{1}{2}}<+\infty ,$$
and the space $ H^{1}_{loc,u}=\{v \mid v,|\nabla v|\in
L^{2}_{loc,u}\}.$
 Moreover, we take
 $f:\er\rightarrow \er $ and $g:\er^{2N+2}\rightarrow \er $  two
${\cal {C}}^1$ functions satisfying  the following conditions:
\begin{eqnarray*}
(H_f)& |{f(u)}|\le M\Big(1+\frac{|u|^{p_c}}{\log^a (2+u^2)}\Big), \
&{\textrm {for all }}\ u\in \er
 \qquad{{\textrm {with}}}\ \ (M>0,\ \  a>1),\\
(H_g)& |{g(x,t,v,z)}|\le M(1+|v|+|z|), &  {\textrm {for all }}\
x,v\in \er^N, t,z\in \er \  \ {{\textrm {with}}}\  (M>0).
\end{eqnarray*}

\medskip

The Cauchy problem of  $\eqref{gen}$ is well posed in
$H^{1}_{loc,u}\times L^{2}_{loc,u}$. This follows from the finite
speed of propagation and the well-posdness in $H^{1}\times L^{2}$,
valid whenever $1< p <p_{S}=1+\frac{4}{N-2}$. 
The existence of
blow-up solutions $u(t)$ of $\eqref{gen}$ follows from
 ODE techniques or the energy-based blow-up criterion by Levine \cite{Ltams74} (see
also \cite{LT},  \cite{T1}).  More
blow-up results can be found in Caffarelli and Friedman
 \cite{CFarma85},  \cite{CFtams86}, Kichenassamy and Littman
 \cite{KL1cpde93}, \cite{KL2cpde93}.
Let us mention the rather surprising result
of Killip and Vi\c san who proved in \cite{KV11} that the ``first''
blow-up set $\{ x_0\;|\;T(x_0) = \min_{x\in \R}T(x)\}$ can be any
Cantor set.
Numerical simulations of blow-up are given by Bizo\'n {\it and al.}  (see \cite{Bjnmp01}, \cite{BBMWnonl10}, \cite{BCTnonl04}, 
\cite{BZnonl09}). 

\medskip

In this paper, we will consider
   $u(x,t)$  a blow-up solution to equation $\eqref{gen}$. We define  $\Gamma =\{(x,T(x))\}$  such that 
the 
  maximal influence domain $D_{u}$  of $u$  is written as
$$D_{u}=\{ (x,t)| t< T(x)\}.$$
Moreover, from the finite speed of propagation, $T$ is a 1-Lipschitz
function. 
The surface $\Gamma$ is called the blow-up graph of $u$. 
A point $x_0\in \R^N$  is a non-characteristic point if there are
\begin{equation}\label{cone}
 \delta_{0}=\delta_{0}(x_{0})\in (0,1)\,
\mbox{ and }t_0<T(x_0)
\  {\rm such }\,\,{\rm that} \,\,u\,\, {\rm is}\,\,
  {\rm defined}\,\, {\rm on}\,\, {\mathcal C}_{x_0, T(x_0), \delta_0}\cap \{t\ge t_0\}
\end{equation}
where 
${\mathcal C}_{\bar x, \bar t, \bar \delta}=\{(x,t)\;|\; t< \bar t-\bar \delta|x-\bar x|\}.$

\medskip

In the case $(f, g) \equiv (0, 0)$,  equation $\eqref{gen}$ reduces
to the semilinear wave equation:
\begin{equation} \label{0}
\partial_{t}^2 u =\Delta u+|u|^{p-1}u,\,\,\,(x,t)\in \R^N \times [0,T).
\end{equation}
It is interesting to recall that   when $1<p\le p_c,$
 in \cite{MZajm03}, \cite{MZimrn05}  and  \cite{MZma05}
Merle and Zaag have proved, that if $u$  is a solution of $\eqref{0}$
with blow-up graph $\Gamma : \{ x\mapsto T(x)\}$ and $x_{0}$ is a
non-characteristic point, then for all $t\in
[\frac{3T(x_{0})}{4},T(x_{0})]$,
 \begin{eqnarray}\label{mzmz}
  0 < \varepsilon_{0}(N,p)\leq (T(x_{0})-t)^{\frac{2}{
   p_c-1}}\frac{\|u(t)\|_{L^{2}(B(x_{0},T(x_{0})-t))}}{(T(x_{0})-t)^{\frac{N}{2}}}\\
  +(T(x_{0})-t)^{\frac{2}{ p_c-1}+1}\Big(\frac{\|\partial_{t} u(t)\|_{L^{2}(B(x_{0},T(x_{0})-t))}}
  {(T(x_{0})-t)^{\frac{N}{2}}}
  +
 \frac{\|\nabla u(t)\|_{L^{2}(B(x_{0},T(x_{0})-t))}}{(T(x_{0})-t)^{\frac{N}{2}}}\Big)\leq K,\nonumber
  \end{eqnarray}
 where the constant $K$ depends only on $N$, $p$ and on an upper bound on
 $T(x_{0})$, ${1}/{T(x_{0})}$, $\delta_{0}(x_{0})$ and the initial data
 in $H^{1}_{loc,u}(\R^N)\times L^{2}_{loc,u}(\R^N)$.

\medskip

 Restricting to the one-dimensional case, Merle and Zaag  fully described the blow-up dynamics for solutions of \eqref{0}
(see 
\cite{MZjfa07}, \cite{MZcmp08}, \cite{MZajm11}, \cite{MZisol10}).
Later, C\^ote and Zaag \cite{CZmulti11} proved the existence of
multi-solitons near characteristic points and gave further
refinements for general solutions of \eqref{0}, still in one space dimension.
Among other results, Merle and Zaag proved that characteristic
points are isolated and that the blow-up set $\{(x,T(x))\}$ is
${\cal {C}}^1$ near non-characteristic points and corner-shaped near
characteristic points. 
In higher dimensions, the method used in  the one-dimensional case 
does not remain valid
 because there  is no   classification of selfsimilar solutions of equation \eqref{gen} in the energy space. However, in the radial case outside the origin, we reduce to  the one-dimensional case with  perturbation and could obtain the same results as for $N=1$ (see \cite{MZbsm11}
and also
the extension by Hamza and Zaag in \cite{HZkg12}
to the Klein-Gordon equation and other
damped lower-order perturbations of equation
 \eqref{0}). Recently, Merle and Zaag could address the higher dimensional case in the subconformal
case and prove the stability of the explicit selfsimilar profile with respect to the blow-up point
and initial data (see 
\cite{MZods}, \cite{MZods14}).
Considering the behavior of radial solutions at the origin, Donninger and Sch{\"o}rkhuber were able to prove the stability of the ODE solution 
$u(t) =
\kappa_0 (p)(T - t)^{-\frac{2}{p-1}} $
 with respect to small perturbations in initial data, in the Sobolev subcritical range \cite{DSdpde12} and also in the supercritical range in \cite{DStams14}. Let us also mention that 
Killip, Stoval and Vi\c san proved in \cite{KSVma14}
that
in  superconformal and Sobolev subcritical range, an upper bound on the blow-up rate is available.
This was further refined by Hamza and Zaag in \cite{HZdcds13}.

\medskip

The question of the perturbed nonlinear wave equation was later
investigated by  Hamza and Zaag in \cite{HZnonl12} and \cite{HZjhde12} 
where they   consider a  class of
 perturbed equations, with $(H_{f})$ and $(H_{g})$ replaced by a
 more restrictive conditions: $|f(u)|\leq M(1+|u|^q)$ and $|g(u)|\leq M(1+|u|)$
  for some $q<p,\,\,M > 0$. Then, they proved a similar result to $\eqref{mzmz}$,
valid in the subconformal case.
  Let us also mention that in  \cite{omar1}, the authors extended the
results obtained in  \cite{HZjhde12} and \cite {HZnonl12}  to 
strong perturbed equation of (\ref{gen})  satisfying $(H_{f})$ and $(H_{g})$
 in the subconformal case ($1<p<p_c$). Recently, in  \cite{omar2}
we extended  the  results known in   \cite{omar1} 
 to the conformal case ($p=p_c$) only for
 radial   solutions, assuming  that the  parameter $a$ satisfies   $a>2$.  However,
 these two assumptions appeared to us as technical and  non
natural.  As a matter of fact, coming with  new ideas  (the use of a Pohozaev identity), we aim 
in this work to remove the radial assumptions, though keeping the condition $a>2$. In fact, 
our main contribution in this paper is to construct a 
  Lyapunov  functional in similarity variables  for the problem \eqref{gen} in the non-radial case, relying on the use of a Pohozaev type identity.

\medskip

Pohozaev type identity has been widely used in mathematics literature and the first results
are due to
\cite{Po1}, where
among other things, he proved the nonexistence
of positive solutions for the elliptic equation $\Delta u+|u|^{p-1}u=0$, in the
supercritical case. Later Giga and Kohn, in \cite{GK2015} characterize all stationary solutions in self-similar variables
of nonlinear heat equations $\partial_t u=\Delta u+|u|^{p-1}u$, in the
subcritical case. Recently the same type of
identity have been used in the analysis of elliptic PDEs
(see \cite{DS1}, \cite{V2015}).
In our work, we construct a Pohozaev identity obtained by multiplying the equation \eqref{A} by $y.\grad w$ in a suitable weighted space. As we see above, the use of this Pohozaev identity is crucial to construct a Lyaponov functional.

\medskip

 Let us introduce the following similarity variables, for any $(x_0,T_0)$ such that $0< T_0\le T(x_0)$:
\begin{equation}\label{scaling}
y=\frac{x-x_0}{T_0-t},\qquad s=-\log (T_0-t),\qquad u(x,t)=(T_0-t)
^{-\frac{2}{p_c-1}}w_{x_0,T_0}(y,s).
\end{equation}
From (\ref{gen}), the  function $w_{x_0,T_0}$  (we write $w$ for
simplicity) satisfies the following equation for all $y\in B\equiv
B(0,1)$ and $s\ge -\log T_0$:
\begin{eqnarray}\label{A}
\partial^2_{s}w &=&div(\nabla w- (y.\nabla w)y)-\frac{2(p_{c}+1)}
{(p_{c}-1)^2}w+|w|^{p_{c}-1}w-\frac{p_{c}+3}{p_{c}-1}\partial_{s}w \\
&&-2y.\nabla
\partial_{s}w+e^{\frac{-2p_{c}s}{p_{c}-1}}f(e^{\frac{2s}{p_{c}-1}}w)+
e^{-\frac{2p_cs}{p_c-1}}g\Big(e^{\frac{(p_c+1)s}{p_c-1}}(\partial_sw+y.\grad w+\frac{2}{p_c-1}w)\Big).\nonumber
\end{eqnarray}
 This change of variables
transforms the backward light cone with vortex $(x_0,T_0)$ into the
infinite cylinder $(y,s)\in B\times [-\log T_0,+\infty)$. In the new
set of variables $(y,s),$ the behavior of $u$ as $t \rightarrow T_0$
is equivalent to the behavior of $w$ as $s \rightarrow +\infty$. In
order to keep our analysis clear, we may assume that
$f(u)=\frac{|u|^{p_c}}{\log^a (2+u^2)}$ and $g \equiv 0$, in the
equation \eqref{gen}
 and refer the reader  to  \cite{HZnonl12} and
\cite{HZjhde12} for straightforward adaptations to
 the general case where $f(u)\not\equiv \frac{|u|^{p_c}}{\log^a (2+u^2)}$ and $g \not\equiv 0$. Also, if $T_0=T(x_0)$, then we simply write $w_{x_0}$
instead of $w_{x_0,T(x_0)}$.

\medskip

 The equation (\ref{A}) will be studied in the Hilbert  space $\cal H$
$${\cal H}=\Big \{(w_1,w_2), |
\displaystyle\int_{B}\Big ( w_2^2 +|\grad w_1|^2(1-|y|^2)+w_1^2\Big) {\mathrm{d}}y<+\infty \Big \}.$$

\medskip

In general,
the treatment of the conformal case  requires a new idea as compared to  
 the subconformal case. In fact, 
the method of perturbation of the Lyapunov functional  used in \cite{HZnonl12} and \cite{omar1} works in the   sub-conformal case  but does not  work in the conformal case. Let us recall that  in \cite{HZjhde12}, we 
studied  the problem   in the conformal case,  if we replace   $(H_{f})$  by a
 more restrictive condition: 
\begin{equation}\label{8jan}
|f(u)|\leq M(1+|u|^q),\quad {\textrm{  for some}} \quad q<p_c.
\end{equation} 
We  proceeded in two steps to construct the  Lyapunov functional: first, we exploited  some functional  to obtain a rough estimate to the blow-up solution namely an 
exponentially large bound. Even though this estimate seems bad, it was very useful to
 allow us to derive  a natural Lyapunov functional for equation \eqref{A}, a  crucial step to derive 
 the optimal estimate as in 
\eqref{mzmz}. Let us note that  
the method used in \cite{HZjhde12} 
under the restrictive condition  \eqref{8jan}
 breaks down when our perturbation is stronger, namely when
$f(u)\equiv \frac{|u|^{p_c}}{\log^a (2+u^2)}$.
Let us mention that  we  overcome this difficulty  with Saidi in \cite{omar2} by proceeding  in three steps
to construct the  Lyapunov functional: first,  as in \cite{HZjhde12} we use  some functional  to obtain an exponentially large estimate for  the blow-up solution. 
Then, we use this exponential bound  to obtain a polynomial estimate. 
However, this step  works only if the solution is radial. Finally this polynomial estimate 
 allows us  to prove that  we have a natural Lyapunov functional for equation \eqref{A}, valid only when $a>2$. Then, we  derive the optimal result \eqref{mzmz} if the solution is radial.  
That obstruction fully justifies our new paper, where we invent a new idea to get our optimal result for a non-radial blow-up solution of $\eqref{A}$, when $a>2$.

\medskip

Let us first recall the rough exponential space-time estimate of the solution $u$   of \eqref{gen} near any
non characteristic point obtained in 
\cite{omar2}. More precisely, we established the
following results:

\medskip

\noindent {\it{\bf { (Exponential space-time estimate of solution of
(\ref{A}).}} Let  $u $   a solution of ({\ref{gen}}) with
blow-up graph $\Gamma:\{x\mapsto T(x)\}$ and  $x_0$ is a non
characteristic point. Then for all $\eta \in (0,1)$,  there exists
$t_{0}(x_{0})\in [ 0,T(x_{0}))$ such that,  for all $s\geq -\log (T(x_{0})-t_{0}(x_{0}))$, we have
\begin{equation}\label{co2}
\int_{s}^{s+1}\int_{B}\Big(|\grad
w_{x_0}(y,\tau)|^2+|w_{x_0}(y,\tau)|^{p_c+1}\Big){\mathrm{d}}y{\mathrm{d}}\tau
\leq K_1 e^{  \eta s},
\end{equation}
and
\begin{equation}\label{co1}
\int_{s}^{s+1}\int_{B}\frac{(\partial_{s}w_{x_0}(y,\tau))^2}{(1-|y|^2)^{1-\eta}}{\mathrm{d}}y{\mathrm{d}}\tau
\leq K_1e^{ \frac{p_c+3}{2} \eta s},
\end{equation}
 where the constant $K_1$ depends only on $N, p, M, a, b, \delta_{0}(x_{0})$, 
 $T(x_{0})$ and \\
$\|(u(t_0(x_0)),\partial_tu(t_0(x_0)))\|_{H^{1}\times
L^{2}(B(x_0,\frac{T(x_0)-t_0(x_0)}{\delta_0(x_0)}) )}$.}

\medskip

In this paper,
by exploiting  a uniform version of the exponential estimates of \eqref{co2} and \eqref{co1} (see
\eqref{co2bis} and \eqref{co1bis}  below), we obtain the following polynomial
space-time
estimate:

\begin{thm}{\bf (A polynomially space-time  estimate of  solution of  \eqref{A})}\label{t1}.\\
Let $u $    a solution of ({\ref{gen}})
with blow-up graph $\Gamma:\{x\mapsto T(x)\}$ and  $x_0$ is a non
characteristic point. Then for all $b\in (1,a)$  there exist  $t_1(x_0)=t_1(x_0,b)\ge t_0(x_0)$ such that
 for all
 $s\ge -\log (T(x_0)-t_1(x_0))$,
\begin{equation}\label{a1}
\int_{s}^{s+1} \!\!\int_{B}\!\!\Big((\partial_s
w_{x_0}(y,\tau))^2+|\grad
w_{x_0}(y,\tau)|^2
+|w_{x_0}(y,\tau)|^{p_c+1}\Big){\mathrm{d}}y{\mathrm{d}}\tau\le
K_2s^b,
\end{equation}
where the constant $K_2$ depends only on $N,  p, M, a, b,  K_1, \delta_{0}(x_{0})$,  $ T(x_{0})$\\
 and 
$\|(u(t_1(x_0)),\partial_tu(t_1(x_0)))\|_{
H^{1}\times
L^{2}(B(x_0,\frac{T(x_0)-t_1(x_0)} {\delta_0(x_0)}) )}$.
\end{thm}
Theorem \ref{t1} can be written in  the original variables in the following corollary:
\begin{cor}\label{09avr}
Let $u$    a solution of ({\ref{gen}})
with blow-up graph $\Gamma:\{x\mapsto T(x)\}$ and  $x_0$ is a non
characteristic point. Then for all $b\in (1,a)$  there exist  $t_1(x_0)=t_1(x_0,b)\ge t_0(x_0)$ such that
 for all
  $t\in [t_1(x_0),T(x_0))$, we have
\begin{equation}\label{t2}
\int_{T(x_0)-t}^{T(x_0)-\frac{t}2}\!\int_{B(x_0,{T(x_0)-\tau })}\!
\!\Big(|\partial_tu(x,\tau )|^2+|\grad u(x,\tau )|^{2}+| u(x,\tau )|^{p_c+1}
\Big)\mathrm{d}x
\mathrm{d}\tau\le K_2 |\log(T(x_0)-t)|^b.\no
\end{equation}
\end{cor}

\medskip

  Now, we are able to adapt the analysis
performed in [9] for equation \eqref{A} and  announce our main result valid only when $a>2$:

\begin{thm}\label{t2}
{\bf {(Blow-up rate for equation \eqref{gen})}}.\\
Let $a>2$, consider   $u $    a solution of ({\ref{gen}}) with blow-up graph
$\Gamma:\{x\mapsto T(x)\}$ and  $x_0$ is a non characteristic point,
then there exist  $\widehat{S}_2$  large enough  such that

i)
 For all
 $s\ge \widehat{s}_2(x_0)=\max(\widehat{S}_2,-\log \frac{T(x_0)}4)$,
\begin{equation*}
0<\varepsilon_0\le \|w_{x_0}(s)\|_{H^{1}(B)}+ \|\partial_s
w_{x_0}(s)\|_{L^{2}(B)} \le K,
\end{equation*}
where $w_{x_0,T(x_0)}$ is defined in (\ref{scaling}) and $B$  is the unit ball of $\er^N$.\\
ii)  For all
  $t\in [t_2(x_0),T(x_0))$, where  $t_2(x_0)=T(x_0)-e^{-\widehat{s}_2(x_0)}$, we have
\begin{eqnarray*}
&&0<\varepsilon_0\le (T(x_0)-t)^{\frac{2}{p_c-1}}\frac{\|u(t)\|_{L^2(B(x_0,{T(x_0)-t}))}}{ (T(x_0)-t)^{\frac{N}{2}}}\nonumber\\
&&+ (T(x_0)-t)^{\frac{2}{p_c-1}+1}\Big
(\frac{\|\partial_tu(t)\|_{L^2(B(x_0,{T(x_0)-t}))}}{
(T(x_0)-t)^{\frac{N}{2}}}+
 \frac{\|\grad u(t)\|_{L^2(B(x_0,{T(x_0)-t}))}}{ (T(x_0)-t)^{\frac{N}{2}}}\Big )\le K,
\end{eqnarray*}
where
$K=K(K_2, N, p, M, a, b, T(x_0), t_2(x_0),\|(u(t_2(x_0)),\partial_tu(t_2(x_0)))\|_{
H^{1}\times
L^{2}(B(x_0,\frac{T(x_0)-t_2(x_0)}{\delta_0(x_0)}) )})$.
\end{thm}


\noindent{\bf Remark 1.1.}\\
Please note that we crucially need a covering technique in our argument, that is why
we need to prove  a uniform version for $x$ near $x_0$ (see   
the exponential space-time estimate  written in \eqref{co2bis}
and \eqref{co1bis},  Theorem 1.1' and  Proposition 3.1.
It happens that the generalization to a uniform version
 valid in  the set
$\{ (x,T_0-\delta_0(x_0)(x-x_0)), t_2(x_0)\le T_0\le T(x_0)\ \  \textrm{and} \ \  |x-x_0|\le \frac{T_0}{\delta_0(x_0)}\}$
is straightforward and  we  refer to \cite{HZjhde12}  for more  details.\\

\medskip

This paper is organized as follows:  In section 2,   we give a new decreasing  functional  for equation \eqref{A}. Then we prove Theorem \ref{t1}, where we obtain 
a  polynomial space-time estimate  of the solution $w$. Using this result,   we prove in section 3 that the "natural" 
functional is a Lyaponov functional for equation \eqref{A} with  the additional  assumption  $a>2$. 
Then, proceeding  as in \cite{omar2},  we prove
Theorem \ref{t2}.\\

We mention that $C$ will be used to denote a constant that's depends on $N$, $a$, $b$ and  $M$ which may vary from line to line. We also introduce
\begin{equation}\label{90}
 F(u)=\int_{0}^{u}f(v){\mathrm{d}}v.
 \end{equation}

\medskip

{\bf Acknowledgment:} 
The author would like  thank the
reviewers for their  valuable comments which undoubtedly helped us to improve the
presentation of our results.


\section{Proof of Theorem \ref{t1}} 
Note that our approach in this section is very close to \cite{omar2}. In fact,
we use  an uniform version for $x$ near $x_0$ for  the exponential bound  on time average to obtain  
an uniform version for $x$ near $x_0$ 
polynomial  estimate on time average of the $H^{1}\times L^{2}(B)$
norm of $(w, \partial_{s} w)$. 
More precisely,
this section is devoted 
to the  proof of  a general version of Theorem \ref{t1}, uniform for $x$ near $x_0$ (see Theorem 1.1'). This
section is divided into four subsections:
\begin{itemize}
\item In the first one we give some classical energy estimates
following  
 from  the multiplication of  equation \eqref{A} by $w(1-|y|^2)^{s^{-b}}$ and $\partial_sw(1-|y|^2)^{s^{-b}}$.
\item The second subsection is devoted to give   new energy estimates  following from the  multiplication of  equation \eqref{A} by $y.\grad w(1-|y|^2)^{1+s^{-b}}(1-\log (1-|y|^2)$.
\item  By combining the above  energy estimates  obtained in the two subsections, we construct a decreasing functional for equation \eqref{A}
and a blow-up criterion involving this functional.
\item Then, we conclude the proof of Theorem 1.1'. 
\end{itemize}
 Now, we start by  stating
the uniform version of the  exponential bound  on time average  for $x$ near $x_0$ obtained in \cite{omar2}.

\medskip

\noindent {\bf{ (Uniform exponential space-time estimate of solution of
(\ref{A}).})} Consider $u $   a solution of ({\ref{gen}}) with
blow-up graph $\Gamma:\{x\mapsto T(x)\}$ and  $x_0$  a non
characteristic point. Then for all $\eta \in (0,1)$,  there exists
$t_{0}(x_{0})\in [ 0,T(x_{0}))$ such that,   for all $T_0 \in [t_{0}(x_0),T(x_{0})]$,  for all $s\geq -\log (T_{0}-t_{0}(x_{0}))$
and $x\in \er^N$, where $|x-x_0|\le \frac{e^{-s}}{\delta_0(x_0)}$,
we have
\begin{equation}\label{co2bis}
\int_{s}^{s+1}\int_{B}\Big(|\grad
w(y,\tau)|^2+|w(y,\tau)|^{p_c+1}\Big){\mathrm{d}}y{\mathrm{d}}\tau
\leq K_1 e^{  \eta s},
\end{equation}
and
\begin{equation}\label{co1bis}
\int_{s}^{s+1}\int_{B}\frac{(\partial_{s}w(y,\tau))^2}{(1-|y|^2)^{1-\eta}}{\mathrm{d}}y{\mathrm{d}}\tau
\leq K_1e^{ \frac{p_c+3}{2} \eta s},
\end{equation}
where $w=w_{x,T^*(x)}$ is defined in \eqref{scaling} with 
 \begin{equation}\label{18dec1}
T^*(x)=T_0-\delta_0(x_0)(x-x_0),
\end{equation}
 $K_1$ depends on $N, p, M, a,  \delta_{0}(x_{0})$, 
 $T(x_{0})$ and 
$\|(u(t_0(x_0)),\partial_tu(t_0(x_0)))\|_{H^{1}\times
L^{2}(B(x_0,\frac{T(x_0)-t_0(x_0)}{\delta_0(x_0)}) )}$.

\bigskip

 Consider $u $   a solution of ({\ref{gen}}) with
blow-up graph $\Gamma:\{x\mapsto T(x)\}$ and  $x_0$ is a non
characteristic point.
Let   $T_0\in (t_0(x_0), T(x_0)]$, for all 
$x\in \er^N$  such that $|x-x_0|\le \frac{T_0}{\delta_0(x_0)}$, then we write $w$ instead of $w_{x,T^*(x)}$ defined in (\ref{scaling}) with $T^*(x)$ given in  (\ref{18dec1}).
As in \cite{omar2},  for any $b\in (1,a)$,  we 
put the equation  in $w$ in the following form:
\begin{eqnarray}\label{var}
\partial_{s}^2w&=&\frac{1}{\phi}\div(\phi \grad w-(y.\grad w)\phi
y)+\frac{2}{s^b}y.\grad w-\frac{2p_c+2}{(p_c-1)^2}w+|w|^{p_c-1}w\\
&&-\frac{p_c+3}{p_c-1}\partial_s w-2y.\grad \partial_sw+e^{\frac{-2p_c s}{p_c-1}}f(e^{\frac{2s}{p_c-1}}w),\
  \forall y\in B\ {\textrm{ and}} \ s\ge -\log T^*(x),\nonumber
\end{eqnarray}
where   
\begin{equation}\label{phi}
\phi=\q=(1-|y|^2)^{s^{-b}}.
\end{equation}

\medskip

 A key step is to find a  funcional  $E(s)$ satisfying a differential inequality
 of type:
 \begin{equation}\label{jardin1}
\frac{d}{ds}E(s)\le -\frac{1}{s^{\mu}}\iint \frac{(\partial_{s}w)^2}{(1-|y|^2)^{1-s^{-b}}}\y
+ \frac{C}{s^{\mu}}E(s), \quad  {\textrm{for some}}\quad  \mu>1.
\end{equation} 

\medskip

In order to control the perturbative terms, we  view the equation \eqref{A}  as a perturbation of the
conformal case  (corresponding to $\phi\equiv 1$ already treated in
\cite{HZjhde12})
  with this term $\frac2{s^b} y \cdot \nabla w$.  Even the term is 
   a lower order term with respect to the nonlinearity this term has a clear effect because an estimate of type
 \eqref{jardin1} implies  
polynomial estimate on time average. 
 It's worth noticing here that the weight $\q$ defined in \eqref{phi}
 (we write $\phi$ for simplicity), depends on time,
 it is not the case in this series of papers
 \cite{MZajm03}, \cite{MZimrn05},
\cite{MZma05}, \cite{MZjfa07}, \cite{MZcmp08}, \cite{MZbsm11}, \cite{HZjhde12}, \cite{HZnonl12}, \cite{HZdcds13}
and \cite{omar1}   we expect that the derivations in time are problematic. In fact,
    we note after observation, that there are new terms appearing compared with the previous works. 
Note that this problem was overcome in the radial case in \cite{omar2},
since the analysis  uses the fact that
 the  tangential
part $\grad_{\theta}w=\grad w-\frac{y.\grad w}{|y|^2}y$ of $\grad w$ vanishes. 
Here, we further refine our argument  allowing  to handle the tangential
part   $\grad_{\theta}w$
 to construct a  function satisfying \eqref{jardin1} in the non-radial case. 
Our method  uses
 a new functional obtained by multiplying the 
equation \eqref{var} by $y.\grad w(1-|y^2)^{1+s^{-b}}\Big(1-\log(1-|y|^2)\Big)$. 
The addition of this functional to 
 some energy estimates 
established by multiplying the equation \eqref{A} by $w$ and $\partial_s w$ in suitable weighted  spaces 
 permitted the control of the bad terms  and is required even in the non-radial case.

\medskip

Notice that in the rest of this section in spirit lightening the paper,  we  define 
\begin{equation}\label{24nov1}
 \nabla_r w=\frac{y.\nabla w}{|y|^2} y\quad \ {\textrm{and}}\quad  \nabla_{\theta}w=
\nabla w-
\frac{y.\nabla w}{|y|^2} y.
\end{equation}
Then, it is given by  \eqref{24nov1}, we can write $\nabla w=\nabla_rw+\nabla_{\theta}w$ and
we have  the identities 
\begin{equation}\label{12nov5}
|y|^2|\grad w|^2-(y.\grad w)^2=|y|^2|\nabla_{\theta}w|^2,
\end{equation}
and
\begin{equation}\label{12nov6}
|\grad w|^2-(y.\grad w)^2=|\nabla_{\theta}w|^2+(1-|y|^2|)|\nabla_{r}w|^2.
\end{equation}

\subsection{Classical  energy  estimates} 
To control the norm of $(w(s), \partial_s w(s))$, 
we start by introducing the following natural 
functionals, for all  $b\in (1,a)$,  
\begin{eqnarray}\label{F0}
E(w(s),s)\!\!\!&=&\!\!\!\!\iint \Big(\frac{1}{2}(\partial_{s}w)^2+\frac{1}{2}(|\grad w|^2-(y.\grad w)^2)+\frac{p_c+1}{(p_c-1)^2}w^2-\frac{|w|^{p_c+1}}{p_c+1}\Big)\q \y\quad\nonumber\\
&&-e^{\frac{-2(p_c+1)s}{p_c-1}} \iint  F(e^{\frac{2s}{p_c-1}}w)\q \y,\no\\
&&\nonumber\\
J(w(s),s)&=&-\frac{1}{s^{b}}\iint   w\partial_{s}w
\q \y+\frac{N}{2s^{b}}\iint   w^2
\q \y, \\
&&\no\\
H(w(s),s)&=&E(w(s),s)+J(w(s),s).\nonumber
\end{eqnarray}

\medskip

 In order to  bound
the time derivative of $H(w(s),s)$, we begin with
 bounding
the time derivative of $E(w(s),s)$  in the following lemma:
\begin{lem} For all  $b\in (1,a)$, $\e\in (0,\frac12)$ and
  $s \geq \max(-\log T^*(x), 1)$, we have 
\begin{eqnarray}\label{jardin2}
\frac{d}{ds}(E(w(s),s))&= &-\frac{2}{s^{b}}\iint (\partial_{s}w)^2\frac{|y|^2\q}{1-|y|^2}\y+\frac{2}{s^{b}}\iint  y.\grad w\partial_{s}w\q \y\\
&& +\frac{b}{s^{b+1}}\iint \Big(\frac{|w|^{p_c+1}}{p_c+1}+e^{-\frac{2(p_c+1)s}{p_c-1}}F(e^{\frac{2s}{p_c-1}}w)\Big)\q \log(1-|y|^2) \y\no\\
&&-\frac{b}{2s^{b+1}}\iint  |\grad_{\theta} w|^2|y|^2\q\log(1-|y|^2) \y
+\Sigma_{1}(s),\no
 \end{eqnarray}
where $\Sigma_{1}(s)$ satisfies
\begin{eqnarray}\label{5novlem}
\Sigma_1(s)  &\le &
\frac{b\e }{2s^{b}}\int_{B}\Big((\partial_{s}w)^2
+|\grad w|^2(1-|y|^2)+\frac{2p_c+2}{(p_c-1)^2} w^2\Big)\q  \y\\
&&+\frac{C }{s^a}\iint |w|^{p_c+1}\q \y
+    C e^{-\frac{\e}{4} s}\int_{B}\Big(\frac{(\partial_{s}w)^2}{\sqrt{1-|y|^2}}
+ |\grad w|^2 +w^2\Big)\y
+Ce^{-s}.\no
\end{eqnarray}
\end{lem}
{\it Proof}: Multiplying $\eqref{var}$ by $\partial_{s} w\!\ \q$ and integrating over the ball $B$, we obtain $\eqref{jardin2}$ where 
 $\Sigma_{1}(s)=\Sigma^1_{1}(s)+\Sigma^2_{1}(s)$ and where
\begin{eqnarray}
\Sigma_{1}^{1}(s)\!\!\!&=&\!\frac{2(p_c+1)}{p_c-1}e^{-\frac{2(p_c+1)s}{p_c-1}}
\iint \!\! \Big( F(e^{\frac{2s}{p_c-1}}w) -\frac{e^{\frac{2s}{p_c-1}}}{p_c+1} 
 wf(e^{\frac{2s}{p_c-1}}w)\Big)\q\y,\ \ \ \quad \label{4nov1}\\
\Sigma_{1}^{2}(s)&=&-\frac{b}{2s^{b+1}}\iint (\partial_{s}w)^2\q\log(1-|y|^2) \y\no\\
&&-\frac{b}{s^{b+1}}\frac{p+1}{(p-1)^2}\iint  w^2\q\log(1-|y|^2) \y\label{4nov2}\\
&&-\frac{b}{2s^{b+1}}\iint |\grad_r w|^2(1-|y|^2)\q \log(1-|y|^2) \y.\no
\end{eqnarray}
Now, we control the terms $\Sigma_{1}^{1}(s)$ and $ \Sigma_{1}^{2}(s)$. 
Clearly the functions $f$ and $F$ defined in \eqref{90} satisfies the following estimate:
\begin{equation}\label{F1nov}
|F(x)|+|xf(x)|\leq C\Big(1+\frac{ |x|^{p_c+1}}{\log^a(2+x^2)}\Big).
\end{equation}
It easily follows from   \eqref{4nov1} and
\eqref{F1nov} that for all $s\geq \max    (-\log T^*(x),1)$, we write
 \begin{equation}\label{4nov3}
 \Sigma_{1}^{1}(s) \leq C\iint \frac{|w|^{p_c+1}}{
\log^a(2+e^{\frac{4s}{p_c-1}}w^2)}
 \q \y+Ce^{-s}.
\end{equation}
As in \cite{omar1} and \cite{omar2}, for all  $s\geq \max  (-\log T^*(x),1)$,
we divide the  ball $B$ into two parts
 \begin{equation}\label{27nov1}
A_{1}(s)=\{y \in B\,\,|\,\, w^2(y,s)\leq  e^{-\frac{2s}{p_{c}-1}}\}\,\,{\rm and }\,\,A_{2}(s)=\{y \in B
\,\,|\,\, w^2(y,s)\ge  e^{-\frac{2s}{p_{c}-1}}\}.
\end{equation}
Using the definition of the set $A_1(s)$ defined in \eqref{27nov1} we get, for all $s\geq \max (-\log T^*(x),1)$
\begin{equation}\label{93}
\int_{A_{1}(s)}\frac{|w|^{p_{c}+1}}{\log^a(2+e^{\frac{4s}{p_{c}-1}}w^2)}\q \y \leq C e^{-\frac{(p_{c}+1)s}{p_{c}-1}}\int_{A_{1}(s)}\q \y\leq C e^{-s}. 
\end{equation}
Also, by using the definition of the set $A_2(s)$ defined in \eqref{27nov1}, we can write   if $y\in A_{2}(s)$, we have $\log(2+e^{\frac{4s}{p_{c}-1}}w^2)\ge \frac{2s}{p_{c}-1},$
 one has,  for all $ s \geq \max (-\log T^*(x),1)$
\begin{equation}\label{102}
\int_{A_2(s)}\frac{|w|^{p_{c}+1}}{\log^a(2+e^{\frac{4s}{p_{c}-1}}w^2)}\q \y 
\leq \frac{C}{s^a}\int_{B}|w|^{p_{c}+1} \q \y.\ \ \
\end{equation}
Hence, the inequality \eqref{4nov3}, \eqref{27nov1}, \eqref{93} and \eqref{102}, imply that
 \begin{equation}\label{pi}
 \Sigma_{1}^{1}(s) \leq \frac{C}{s^a}\iint |w|^{p_c+1} \q \y+ Ce^{-s},\quad  \ {\textrm{for\ all}}\ s \geq \max (-\log T^*(x),1).
\end{equation}
Now, we are going to estimate $\Sigma_{1}^{2}(s)$. 
The treatment of this term is more difficult
because it contains terms  
 with singular weight.  In fact,  unlike the terms in  $\Sigma_{1}^{1}(s)$ which 
constituting
by therms with
the weight is $\q$, here the weight is $-\q \log(1-|y|^2)$.
 To overcome this problem, we divide the unit  ball $B$ in two parts:
a first part which is near the boundary $\partial B$  and the other is  the rest of the unit ball.  More precisely,
let $\e \in (0,\frac12)$, for all 
$s\ge \max  (-\log T^*(x),1)$, we divide $B$ into two parts
 \begin{equation}\label{B1}
 B_{1}(s)=\{y \in B\,\,|\,\, 1-|y|^2\leq  e^{-\e s}\}\,\,{\rm and }\,\,B_{2}(s)=\{y \in B\,\,|\,\, 1-|y|^2\ge  e^{-\e s}\}.
\end{equation}
We see that:
$\Sigma_{1}^{2}(s)=\underbrace{\chi^1_{1}(s)+\chi^2_{1}(s)+\chi^3_{1}(s)}_{\chi_{1}(s)}+\chi_{2}(s),$
where 
\begin{eqnarray*}
\chi^1_{1}(s)&=&-\frac{b}{2s^{b+1}}\int_{B_1(s)}(\partial_{s}w)^2\q\log(1-|y|^2) \y,\\
\chi^2_{1}(s)&=&-\frac{b}{2s^{b+1}}\int_{B_1(s)}|\grad_r w|^2(1-|y|^2)\q\log(1-|y|^2) \y,\\
\chi^3_{1}(s)&=&-\frac{b}{s^{b+1}}\frac{p_c+1}{(p_c-1)^2}\int_{B_1(s)} w^2\q\log(1-|y|^2) \y,
\end{eqnarray*}
and where
\begin{eqnarray*}
\chi_{2}(s)&=&-\frac{b}{2s^{b+1}}\int_{B_2(s)}\Big((\partial_{s}w)^2
+\frac{2p_c+2}{(p_c-1)^2} w^2\Big)\q
\log(1-|y|^2) \y\\
&&-\frac{b}{2s^{b+1}}\int_{B_2(s)}|\grad_r w|^2(1-|y|^2)\q \log(1-|y|^2) \y.
\end{eqnarray*}
From   the fact that, for all $s\ge \max  (-\log T^*(x),1)$,   the function $y\mapsto \q (1-|y|^2)^{\frac{1}{4}}\log(1-|y|^2)$ is   uniformly bounded on $B$  and using  the inequality 
 $$(1-|y|^2)^{\frac{1}{4}}\leq  e^{-\frac{\e}{4} s},\quad \quad 
{\textrm {for all}}\quad y \in B_{1}(s),$$
 we can write, for all $s\ge \max  (-\log T^*(x),1)$,
\begin{equation}\label{X11}
\chi^1_{1}(s)  \leq  C e^{-\frac{\e}{4} s}\int_{B}\frac{(\partial_{s}w)^2}{\sqrt{1-|y|^2}}\y.
\end{equation}
Similarly  we obtain easily, for all $s\ge \max  (-\log T^*(x),1)$,
\begin{equation}\label{X12}
\chi^2_{1}(s) +\chi^3_{1}(s)  \leq  C e^{-\frac{\e}{4} s}\int_{B}|\grad w|^2 \y+
  C e^{-\frac{\e}{4} s}\int_{B}\frac{w^2}{\sqrt{1-|y|^2}}\y.
\end{equation}
 Let us recall from   \cite{MZajm03} the following Hardy type inequality
\begin{equation}\label{Hardy7nov}
 \int_{B}w^2\frac{|y|^2}{\sqrt{1-|y|^2}}\y\leq C\int_{B}|\nabla w|^2(1-|y|^2)^{\frac32} \y+C\int_{B} w^2\sqrt{1-|y|^2} \y.
\end{equation}
Using the fact that $\frac{w^2}{\sqrt{1-|y|^2}}=\frac{|y|^2w^2}{\sqrt{1-|y|^2}}+w^2\sqrt{1-|y|^2}$,
we conclude that 
\begin{equation}\label{X13}
\chi^2_{1}(s) +\chi^3_{1}(s)  \leq  C e^{-\frac{\e}{4} s}\int_{B}w^2\y+C e^{-\frac{\e}4 s}\int_{B}|\grad w|^2 \y.
\end{equation}
Combining \eqref{X11} and  \eqref{X13}, one easily obtain 
\begin{equation}\label{X1}
\chi_1(s)  \leq  C e^{-\frac{\e}{4} s}\int_{B}\Big(\frac{(\partial_{s}w)^2}{\sqrt{1-|y|^2}}
+  |\grad w|^2 
+ w^2\Big)\y.
\end{equation}
Next,   by \eqref{12nov6} and by exploiting the fact that
if $y \in B_{2}(s)$, we have $-\log(1-|y|^2)\leq \e s$, we have
  for all $s\ge \max  (-\log T^*(x),1)$, 
\begin{equation}\label{5nov11}
\chi_{2}(s)\le\frac{b\e }{2s^{b}}\int_{B}\Big((\partial_{s}w)^2
+|\grad w|^2(1-|y|^2)+\frac{2p_c+2}{(p_c-1)^2} w^2\Big)\q  \y.
\end{equation}
Then we infer from  \eqref{X1}, \eqref{5nov11} and the identity $\Sigma^2_1(s)=  \chi_{1}(s)+\chi_{2}(s)$
  we have,  for all $s\ge \max  (-\log T^*(x),1)$
\begin{eqnarray}\label{X5nov}
\Sigma^2_1(s)  &\leq&  
\frac{b\e }{2s^{b}}\int_{B}\Big((\partial_{s}w)^2
+|\grad w|^2(1-|y|^2)+\frac{2p_c+2}{(p_c-1)^2} w^2\Big)\q  \y
\\
&&+  C e^{-\frac{\e}{4} s}\int_{B}\Big(\frac{(\partial_{s}w)^2}{\sqrt{1-|y|^2}}+|\grad w|^2 +w^2\Big)\y.\no
\end{eqnarray}
The result \eqref{5novlem} derives immediately from  \eqref{pi}, \eqref{X5nov} and  the identity $\Sigma_1(s)  =\Sigma^1_1(s)  +\Sigma^2_1(s)$,
which ends the proof of Lemma 2.1.
\Box

 \medskip

We are going to prove the following estimate to the functional $J(w(s),s)$. 
 \begin{lem}
 For all  $b\in (1,a)$ and $\e\in (0,\frac12)$,
there exists $S_{1}\geq 1$ such that for all $s \geq \max  (-\log T^*(x),S_1)$, we have the following inequality:
\begin{eqnarray}\label{6nov1}
\frac{d}{ds}(J(w(s),s))&=&-\frac{p_c+7}{4s^{b}}\iint   (\partial_{s}w)^2\q \y
+\frac{p_c+3}{2s^{b}}H(w(s),s)\\
&&-\frac{p_c-1}{4s^{b}}\iint (|\grad w|^2-(y.\grad w)^2)\q\y
-\frac2{s^b} \iint  \partial_swy.\grad w\q\y\no\\
&&
-\frac{(p_c+1)}{2(p_c-1)s^b} \iint  w^2\q\y
-\frac{p_c-1}{2(p_c+1)s^b} \iint  |w|^{p_c+1}\q\y+\Sigma_2(s),\no
\end{eqnarray}
where $\Sigma_2(s)$ satisfies
\begin{eqnarray}\label{7nov1}
  \Sigma_{2}(s)&\leq& \frac{C}{s^{b+1}}\iint (\partial_{s}w)^2\q\y + \frac{34}{(p_c+16)s^{b}}\iint (\partial_{s}w)^2\frac{|y|^2\q}{1-|y|^2}\y\\
&&+\frac{p_c-1}{8s^{b}}\iint  |\nabla w|^2(1-|y|^2)\q\y
  +\frac{3(p_c+1)}{8(p_c-1)s^b}\iint  w^2\q \y
\no\\
&&+\frac{C}{s^{a+b}} \iint |w|^{p_c+1}\q\y+
 C e^{-\frac{\e}{4} s}\int_{B}\Big(w^2+|\grad w|^2\Big) \y+ Ce^{-s}.\no
\end{eqnarray}
\end{lem}
{\it Proof}: Note that $J(w(s),s)$ is a differentiable function and  we get,  for all $s\geq  -\log T^*(x)$ 
\begin{eqnarray}\label{j1}
\frac{d}{ds}(J(w(s),s))&=&-\frac{1}{s^{b}}\iint   
(\partial_{s}w)^2\q \y-\frac{1}{s^{b}}\iint   w\partial^2_{s}w\q \y\no\\
&& +\big(\frac{b}{s}+N\big)\frac{1}{s^{b}}\iint   w\partial_{s}w\q \y+\frac{b}{s^{2b+1}}
\iint  w\partial_{s}w  \log(1-|y|^2)\q \y\no\\
&&-\frac{Nb}{2s^{b+1}}\iint   w^2
\q \y-\frac{Nb}{2s^{2b+1}}\iint   w^2
\q\log(1-|y|^2) \y.
\end{eqnarray}
 According to equation $\eqref{var}$, we obtain
\begin{eqnarray}\label{j2}
\frac{d}{ds}(J(w(s),s))&=&
-\frac2{s^b} \iint  \partial_swy.\grad w\q\y
-\b \iint  |w|^{p_c+1}\q\y\\
&&-\frac{1}{s^{b}}\iint   
(\partial_{s}w)^2\q \y+
\frac{1}{s^{b}}\iint (|\grad w|^2-(y.\grad w)^2)\q\y\no\\
&&+\frac{2p_c+2}{(p_c-1)^2}
 \b \iint  w^2\q\y
+\Sigma^1_2(s)+\Sigma_{2}^{2}(s)+\Sigma_{2}^{3}(s)+\Sigma_{2}^{4}(s)
+\Sigma_{2}^{5}(s),\no
 \end{eqnarray}
 where
 \begin{eqnarray*}
\Sigma_{2}^{1}(s)&=&
\frac{N(2-b)}{2s^{b+1}}\iint   w^2\q \y+\frac{b}{s^{b+1}} \iint  w\partial_s w\q\y,\no\\
  \Sigma_{2}^{2}(s)&=&-\frac{2}{s^{3b}}\iint  w^2\frac{|y|^2\q}{1-|y|^2}\y
+\frac{4}{s^{2b}} \iint  w\partial_{s}w\frac{|y|^2\q}{1-|y|^2}\y, \\ 
 \Sigma_{2}^{3}(s)&=&-\frac{e^{\frac{-2p_cs}{p_c-1}}}{s^{b}}\iint  wf(e^{\frac{2s}{p_c-1}}w)\q\y,\no\\
\Sigma_{2}^{4}(s)&=&\frac{b}{s^{2b+1}}\iint  w\partial_{s}w \q \log(1-|y|^2) \y,\no\\
\Sigma_{2}^{5}(s)&=&
-\frac{Nb}{2s^{2b+1}}\iint   w^2\q\log(1-|y|^2) \y.
 \end{eqnarray*}
According to the expression of $H(w(s),s)$, with some straighforward computation we obtain \eqref{6nov1} where
 \begin{equation}\label{7nov2}
\Sigma_2(s)=\Sigma^1_2(s)+\Sigma_{2}^{2}(s)+\Sigma_{2}^{3}(s)+\Sigma_{2}^{4}(s)
+\Sigma_{2}^{5}(s)+\Sigma_{2}^{6}(s)++\Sigma_{2}^{7}(s),
\end{equation} 
and where
\begin{eqnarray*}
\Sigma_{2}^{6}(s)&=&\frac{p_c+3}{2s^b} e^{\frac{-2(p_c+1)s}{p-1}} \iint  F(e^{\frac{2s}{p-1}}w)\q \y,\no\\
 \Sigma_{2}^{7}(s)&=&\frac{p_c+3}{2s^{2b}} \iint  w\partial_s w\q\y
-\frac{N(p_c+1)}{4s^{2b}}  \iint  w^2\q\y.\no
 \end{eqnarray*}
We are going now to estimate  the  different terms of \eqref{7nov2}.
Thanks to the  the classical inequality $ab\le a^2+b^2,$  we conclude that  for all $s\ge \max  (-\log T^*(x),1)$
  \begin{equation}\label{R1}
  \Sigma_{2}^{1}(s)+ \Sigma_{2}^{7}(s)\leq \frac{C}{s^{b+1}}\iint \big((\partial_{s}w)^2 
  +  w^2\big)\q \y.
  \end{equation}
By the Cauchy-Schwarz inequality, we write for all $\mu\in (0,1)$ 
\begin{equation}\label{N1}
\Sigma_{2}^{2}(s)\leq \frac{2(1-\mu)}{s^{b}}\iint (\partial_{s}w)^2\frac{|y|^2\q}{1-|y|^2}\y
+\frac{2\mu}{(1-\mu)s^{3b}}\iint  w^2\frac{|y|^2\q}{1-|y|^2}\y.
\end{equation}
Let us recall 
 from   \cite{MZajm03} the following Hardy type inequality, for all $\eta\in (0,1)$
\begin{equation}\label{Hardybis}
 \int_{B}h^2\frac{|y|^2\rho_{\eta}}{1-|y|^2}\y\leq \frac{1}{\eta^2}\int_{B}|\nabla h|^2(1-|y|^2)\rho_{\eta} \y+\frac{N}{\eta}\int_{B} h^2\rho_{\eta} \y.
\end{equation}
Then, from \eqref{Hardybis}, it follows that
\begin{equation}\label{N2}
\iint  w^2\frac{|y|^2\q }{1-|y|^2}\y \leq s^{2b}\iint |\nabla w|^2(1-|y|^2)\q \y+Ns^{b}\iint  w^2\q \y.
\end{equation}
From $\eqref{N1}$, $\eqref{N2}$ with $\mu= \frac{p_c-1}{p_c+15}$, we conclude that,
for all $s\ge \max(-\log T^*(x),1)$
\begin{eqnarray}\label{N3}
\Sigma_{2}^{2}(s)&\leq& \frac{32}{(p_c+15)s^{b}}\iint (\partial_{s}w)^2\frac{|y|^2\q}{1-|y|^2}\y
+\frac{p_c-1}{8s^{b}}\iint  |\nabla w|^2(1-|y|^2)\q\y\no\\
&&+\frac{C}{s^{2b}}\iint  w^2\q\y.
\end{eqnarray}
The same type of estimates used to obtain \eqref{pi} are used here to deduce  easily, for all $s\ge \max  (-\log T^*(x),1)$
 \begin{equation}\label{N4}
\Sigma_{2}^{3}(s)+\Sigma_{2}^{6}(s)\leq \frac{C}{s^{a+b}} \iint |w|^{p_c+1}\q\y+Ce^{-s}.
\end{equation}
Furthermore by using the inequality $ab\le a^2+b^2$  we  write, for all $s\geq \max  (-\log T^*(x),1)$ 
 \begin{equation*}
\Sigma_{2}^{4}(s)\le
\frac{(p_c+1)}{4(p_c-1)s^b}
  \iint  w^2\q\y+
\frac{C}{s^{3b+2}}\iint  (\partial_{s}w)^2  (\log(1-|y|^2))^2\q \y.
\end{equation*}
Using  the fact that,    the function $y\mapsto {|y|^{-2}} \sqrt{1-|y|^2}(\log(1-|y|^2))^2$ is 
 bounded on $B$,  we obtain   
\begin{equation}\label{N5}
\Sigma_{2}^{4}(s)\le 
\frac{(p_c+1)}{4(p_c-1)s^b}
  \iint  w^2\q\y+\frac{C}{s^{3b+2}}\iint   (\partial_{s}w)^2
\frac{|y|^2\q}{\sqrt{1-|y|^2}} \y.
  \end{equation}
Finally,   for all 
$s\ge \max(-\log T^*(x),1)$,
we are going to estimate $\Sigma_{2}^{5}(s)$. For this,  we divide $B$ into two parts $B_{1}(s)$ and $B_{2}(s)$  as defined in \eqref{B1}.
We write
$\Sigma_{2}^{5}(s)=\chi_{3}(s)+\chi_{4}(s),$
where 
\begin{eqnarray*}
\chi_3(s)&=&-\frac{Nb}{2s^{2b+1}}\int_{B_1(s)}   w^2\q\log(1-|y|^2) \y,\no\\
\chi_4(s)&=&-\frac{Nb}{2s^{2b+1}}\int_{B_2(s)}   w^2\q\log(1-|y|^2) \y.\no
\end{eqnarray*}
Hence, if $y \in B_{1}(s)$, then $(1-|y|^2)^{\frac{1}{4}}\leq  e^{-\frac{\e}{4} s}$, 
taking  into account the fact that, for all $s\ge \max (s_0,1)$,   the function $y\mapsto |y|^{-2}\q (1-|y|^2)^{\frac{1}{4}}\log(1-|y|^2)$ is uniformly  bounded on $B$, we can write for all $s\ge \max  (-\log T^*(x),1)$
\begin{equation}\label{X511}
\chi_3(s)  \leq  C e^{-\frac{\e}{4} s}\int_{B}w^2\frac{|y|^2}{\sqrt{1-|y|^2}}\y.
\end{equation}
Furthermore, thanks to   the Hardy-Sobolev inequality \eqref{Hardy7nov},
we conclude that 
\begin{equation}\label{7nov9}
\chi_3(s)  \leq  C e^{-\frac{\e}{4} s}\int_{B}\big(w^2+|\grad w|^2\big) \y.
\end{equation}
Taking in consideration the fact that, if $y \in B_{2}(s)$, we have $-\log(1-|y|^2)\leq \e s$, we obtain
\begin{equation}\label{7nov10}
\chi_{4}(s)\le  
\frac{C}{s^{2b}}\iint   w^2\q\y.
\end{equation}
By adding \eqref{7nov9} and \eqref{7nov10}, we write
\begin{equation}\label{7nov11}
\Sigma_{2}^{5}(s)  \le  \frac{C}{s^{2b}}\iint   w^2\q\y+
 C e^{-\frac{\e}{4} s}\int_{B}\big(w^2+|\grad w|^2\big) \y.
\end{equation}
Consequently, collecting   \eqref{R1}, \eqref{N3},  \eqref{N4},  \eqref{N5}  and  \eqref{7nov11}, one easily
there exists $S_1\ge 1$ such that we obtain $\Sigma_{2}(s)$ satisfies \eqref{7nov1},
 which end the proof of Lemma 2.2.

 \Box

Lemmas 2.1 and 2.2 allows to prove the following lemma:
\begin{lem}
 For all  $b\in (1,a)$ and $0<\e\le   \frac{p_c-1}{32b(p_c+1)}$,
there exist $S_{2}\geq S_1$  and $\lambda_0>0$ such that for all $s \geq \max  (-\log T^*(x), S_{2})$, we have the following inequality:
\begin{eqnarray}\label{7nov33}
\frac{d}{ds}(H(w(s),s))&\le &
\frac{p_c+3}{2s^{b}}H(w(s),s)
{-\frac{b}{2s^{b+1}}\iint  |\grad_{\theta} w|^2|y|^2\q\log(1-|y|^2) \y}
\no\\
&& +\frac{b}{s^{b+1}}
\iint \Big(\frac{|w|^{p_c+1}}{p_c+1}+e^{-\frac{2(p_c+1)s}{p_c-1}}F(e^{\frac{2s}{p_c-1}}w)\Big)\q\log(1-|y|^2) \y\no\\
&&
- \frac{\lambda_0}{s^b}\iint (|\grad w|^2-(y.\grad w)^2)\q\y\\
&&- \frac{\lambda_0}{s^b}\iint \Big(\frac{(\partial_{s}w)^2}{1-|y|^2} 
+ |w|^{p_c+1} +  w^2\Big)\q \y+\Sigma_3 (s),\no
 \end{eqnarray}
where $\Sigma_3(s)$ satisfies
\begin{equation}\label{7nov15}
\Sigma_3(s)  \leq C e^{-\frac{\e}{4} s}\int_{B}\Big(\frac{(\partial_{s}w)^2}{\sqrt{1-|y|^2}}  
+|\grad w|^2 +w^2\Big) \y+ Ce^{-s}.
\end{equation}
\end{lem}

\noindent{\bf Remark 2.1.}
 We notice that  the term $-\frac{b}{2s^{b+1}}\iint  |\grad_{\theta} w|^2|y|^2\q\log(1-|y|^2) \y $ in the inequality (\ref{7nov33})
is non negative, which does not allow to construct 
 a decreasing functional for equation \eqref{var}. 
Please note that in \cite{omar2}    we only treat  the    radial solutions where     this term vanishes. 
One main reason for this restriction  is that
we did not know  control this term in the case of non-radial solutions. 
Here, let us recall that we consider the non-radial case that's why we need  a new idea to overcome this problem. In fact, 
we construct a new functional
$L(w(s),s)$ which is   crucial  to   obtain
 a decreasing functional for equation \eqref{var} later.


\subsection{New  energy  estimates} 
In this subsection, we start by  introduce the   crucial   new
  functional $\L$ defined by the following:
\begin{equation}\label{L0}
L(w(s),s)=\int_{ B}\Big((y.\nabla w)^2+y.\grad w \partial_sw
\Big)\p {\mathrm{d}}y,
\end{equation}
where $\p=(1-|y|^2)\q \big(1-\log(1-|y|^2)\big)$.
As one can see in the statement below, this quantity arises from a Pohozaev identity obtained through 
the multiplication of equation \eqref{A} by $y.\grad w$. This is the main novelty of our paper. More precisely,
to estimate    the time derivative of the
 functional $\L$, we claim   the following:
\begin{lem}\label{L1}
For   all $b\in (1,a)$ and  $s\geq \max  (-\log T^*(x),1)$, we have
\begin{eqnarray}\label{14nov3}
\frac{d}{ds}\L
 &=&(1+\b)\int_{B}|\grad_{\theta} w|^2|y|^2\q\log(1-|y|^2)
{\mathrm{d}}y\\
&&\!\!\!-(2+\frac2{s^b})\!\int_{ B}\! \Big(\frac{|w|^{p_c+1}}{p_c+1}+
e^{\frac{-2(p_c+1)s}{p_c-1}} F(e^{\frac{2s}{p_c-1}}w)\Big) \q \log(1-|y|^2)
\y+\Sigma_3(s),\no
\end{eqnarray}
where
\begin{eqnarray}
\label{cor17nov}
\Sigma_3(s)&\le&
C \int_{ B}  \Big(\frac{(\partial_sw)^2} {1-|y|^2}+|\grad w|^2-(y.\grad  w)^2
+w^2+ |w|^{p_c+1}\Big) \q
\y\\
 && + C e^{-\frac{s}4}\int_{B}|\grad w|^2\y+Ce^{-s}.\no
\end{eqnarray}
\end{lem}
{\it Proof:} Note that  $\L$ is a differentiable
function and we get,  for all
 $s\ge    -\log T^*(x)$  
\begin{equation}\label{20}
\frac{d}{ds}\L=  \int_{B}y.\grad  w (\partial^2_sw+2y.\grad \partial_s w)\p \y+
\Sigma_3^1(s)+\Sigma_3^2(s),
\end{equation}
where
\begin{eqnarray}\label{12nov1}
\Sigma_3^1(s)&=&  \int_{ B}y.\grad \partial_sw \partial_sw \p \y,\no\\
\Sigma_3^2(s)&=&-\frac{b}{s^{b+1}}\int_{ B}\Big((y.\nabla w)^2+y.\grad w \partial_sw
\Big)\p \log(1-|y|^2)\y.\no
\end{eqnarray}
By using \eqref{A} and   integrating by parts, we have
\begin{eqnarray}\label{122}
\frac{d}{ds}\L&=&  \int_{B}y.\grad  w \div(\nabla w-
(y.\nabla w)y)\p\y\\
&&-\int_{ B}\Big(\frac{ |w|^{p_c+1}}{p_c+1}+
e^{\frac{-2(p_c+1)s}{p_c-1}} F(e^{\frac{2s}{p_c-1}}w)\Big) \div(\p y) \y\no\\
&&+\Sigma_3^1(s)+\Sigma_3^2(s)+\Sigma_3^3(s),\no
\end{eqnarray}
where
\begin{equation}\label{14nov44}
\Sigma_3^3(s)= \frac{p_c+1}{(p_c-1)^2}\int_{B}w^2\div(\p y)\y-\frac{p_c+3}{p_c-1}
\int_{B}y.\grad w \partial_{s}w\p \y.
\end{equation}
A straightforward computation yields
 the identity
\begin{equation}\label{12nov3}
\div(\p y)=(N+2+\frac2{s^b})\p -(2+\frac2{s^b})\frac{\p}{1-|y|^2}
+2|y|^2\q.
\end{equation}
From \eqref{12nov3},
we obtain
\begin{eqnarray}\label{14nov1}
-\int_{ B}\Big(\frac{ |w|^{p_c+1}}{p_c+1}+
e^{\frac{-2(p_c+1)s}{p_c-1}} F(e^{\frac{2s}{p_c-1}}w)\Big) \div(\p y) \y\qquad\qquad\qquad\qquad\qquad\qquad \\
=-(2+\frac2{s^b})\int_{ B} \Big(\frac{|w|^{p_c+1}}{p_c+1}+
e^{\frac{-2(p_c+1)s}{p_c-1}} F(e^{\frac{2s}{p_c-1}}w)\Big) \q \log(1-|y|^2)
\y+\Sigma_3^4(s),\no
\end{eqnarray}
where
\begin{eqnarray}\label{14nov2}
\Sigma_3^4(s)&=&-(N
+2+\frac2{s^b})
\int_{ B}\Big(\frac{ |w|^{p_c+1}}{p_c+1}+
e^{\frac{-2(p_c+1)s}{p_c-1}} F(e^{\frac{2s}{p_c-1}}w)\Big) \p  \y\qquad\\
&&+2\int_{ B}\Big(\frac{ |w|^{p_c+1}}{p_c+1}+
e^{\frac{-2(p_c+1)s}{p_c-1}} F(e^{\frac{2s}{p_c-1}}w)\Big)(1+\frac1{s^b}- |y|^2)\q  \y.\no
\end{eqnarray}
After some simple  integration by parts that we leave to appendix A,  \eqref{122} and \eqref{14nov1}, 
we obtain \eqref{14nov3}
where 
\begin{equation}\label{14nov5}
\Sigma_3(s)=\Sigma_3^1(s)+\Sigma_3^2(s)+\Sigma_3^3(s)+\Sigma_3^4(s)+\Sigma_3^5(s),
\end{equation}
and where
\begin{eqnarray}
\Sigma_3^5(s)
 &=&-\b\int_{B}|\grad_{\theta} w|^2|y|^2\q
{\mathrm{d}}y
+\frac{N-2}2\int_{B}\big(|\grad w|^2-(y.\grad  w)^2\big) \p
{\mathrm{d}}y 
\no\\
&&+\b  \int_{B}(y.\grad  w)^2
\p {\mathrm{d}}y- \int_{B}(1-|y|^2)(y.\grad  w)^2
\q {\mathrm{d}}y.
\end{eqnarray}
Now, we    control  all  the   terms on the right-hand side of the identity  
 \eqref{14nov5}.
 After  integration by parts, we use 
 \eqref{12nov3} to show
\begin{equation*}
\Sigma_3^1(s)=-\frac{N}2\int_{ B}  (\partial_sw)^2 \p \y
+(1+\frac1{s^b}) \int_{ B}  (\partial_sw)^2 \frac{|y|^2\p}{1-|y|^2}\y- \int_{ B}  |y|^2 (\partial_sw)^2\q\y.
\end{equation*}
To estimate $\Sigma_3^1(s)$, we using the fact that    $0\le \p\le C\q$, for all  $y\in B$, to write for all 
$s\ge \max  (-\log T^*(x),1)$
\begin{equation}\label{14nova1}
\Sigma_3^1(s)\le
C \int_{ B}  \frac{(\partial_sw)^2}{1-|y|^2}\q\y.
\end{equation}
Note that by  using the inequality   $0\le -\p\log(1-|y|^2)\le C\q$, for all  $y\in B$
and the fact that  $ab\le a^2+b^2$,
 we write for all $s\ge \max  (-\log T^*(x),1)$
\begin{equation}\label{5dec5}
\Sigma_3^2(s)\le \underbrace{-\frac{C}{s^{b+1}}\int_{ B}|\grad w|^2\p \log(1-|y|^2)\y}_{\chi_5(s)}
+\frac{C}{s^{b+1}}\int_{ B}(\partial_sw
)^2\q \y
\end{equation}
We would like now to find an estimate from  the term $\chi_5(s)$. For this, we divide $B$ into two parts 
$B_3(s)$ and $B_4(s)$   defined  for all 
$s\ge \max  (-\log T^*(x),1)$ by
\begin{equation}\label{B431}
 B_{3}(s)=\{y \in B\,\,|\,\, 1-|y|^2\leq  e^{ -s}\}\,\,{\rm and }\,\,B_{4}(s)=\{y \in B\,\,|\,\, 1-|y|^2\ge  e^{- s}\}.
\end{equation}
We write $\chi_5(s)=\chi_5^1(s)+\chi_5^2(s)$, where
\begin{eqnarray*}
\chi^1_{5}(s)&=&-\frac{C}{s^{b+1}}\int_{B_3(s)}|\grad w|^2\p \log(1-|y|^2)\y,\\
\chi^2_{5}(s)&=&-\frac{C}{s^{b+1}}\int_{B_4(s)}|\grad w|^2\p \log(1-|y|^2)\y.
\end{eqnarray*}
From   the fact that,  for all $s\ge \max   (-\log T^*(x),1)$   the function $y\mapsto  \q\sqrt{1-|y|^2}(1-\log(1-|y|^2))
\log (1-|y|^2)$ is uniformly   bounded on $B$  and using  the inequality 
 $$\sqrt{1-|y|^2}\leq  e^{-\frac{s}2},\quad \quad 
\forall y \in B_{3}(s),$$
 we can write, for all $s\ge \max  (-\log T^*(x),1)$,
\begin{equation}\label{5dec1}
\chi^1_{5}(s)  \leq  C e^{-\frac{s}{4} }\int_{B}|\grad w|^2\y.
\end{equation}
Next,  if $y \in B_{4}(s)$, we can write $-(1-\log(1-|y|^2))\log(1-|y|^2)\leq  C s^2$, then we have
  for all $s\ge \max  (-\log T^*(x),1)$, 
\begin{equation}\label{5dec2}
\chi_{5}^2(s)\le\frac{C }{s^{b-1}}\int_{B}
|\grad w|^2(1-|y|^2)\q  \y.
\end{equation}
 By using \eqref{5dec1} and \eqref{5dec2}, we can write for all $s\ge \max  (-\log T^*(x),1)$
\begin{equation}\label{5dec1}
\chi_{5}(s)  \leq  C e^{-\frac{s}{4} }\int_{B}|\grad w|^2\y.
+\frac{C }{s^{b-1}}\int_{B}
|\grad w|^2(1-|y|^2)\q  \y.
\end{equation}
By combining \eqref{5dec5} and \eqref{5dec1} we conclude, for all $s\ge \max  (-\log T^*(x),1)$, 
\begin{equation}\label{5dec6}
\Sigma_3^2(s)\le 
\frac{C}{s^{b-1}}\int_{ B}\big(|\grad w|^2(1-|y|^2)+(\partial_sw)^2\big)\q  \y+C e^{-\frac{s}{4} }\int_{B}|\grad w|^2\y.
\end{equation}
Adding  \eqref{14nov44}  to the identity
\eqref{12nov3},  we can write
\begin{eqnarray}\label{15nov1}
\Sigma_3^3(s)&=&N \frac{(p_c+1)}{(p_c-1)^2}\int_{B}w^2\p\y
-2(1+\frac1{s^b})\frac{(p_c+1)}{(p_c-1)^2}\int_{B}w^2
\frac{|y|^2\p}{1-|y|^2}\y\no\\
&&+2\frac{(p_c+1)}{(p_c-1)^2}\int_{B}|y|^2w^2\q\y
\underbrace{-\frac{p_c+3}{p_c-1}
\int_{B}y.\grad w \partial_{s}w\p \y}_{\chi_6(s)}.
\end{eqnarray}
Note that, by using  inequality $ab\le a^2+b^2$ we obtain
\begin{equation}\label{15nov2}
\chi_6(s) \le \int_{B}|\grad w|^2(1-|y|^2)\p \big(1-\log(1-|y|^2)\big)\y  +C\int_{B}(\partial_{s}w)^2\q \y.
\end{equation}
From \eqref{15nov1}, \eqref{15nov2}, the fact that, for all $y\in B$, we have  $\p\le C \q$ and
 for all $s\ge \max  (-\log T^*(x),1)$ 
   the function $y\mapsto  (1-|y|^2) (1-\log(1-|y|^2))^2$  is uniformly bounded on $B$, we write 
\begin{equation}\label{17nov3}
\Sigma_3^3(s)\le C\int_{B}\big( (\partial_{s}w)^2+w^2
+|\grad w|^2(1-|y|^2)\big)\q \y.
\end{equation}
The same type of estimates used to obtain \eqref{pi} are uesd here together with  the inequality $\p \le C \q$  to deduce  easily, for all $s\ge \max  (-\log T^*(x),1)$ 
\begin{equation}\label{17nov4}
\Sigma_3^4(s)\le
C\int_{ B} |w|^{p_c+1} \q
\y+Ce^{-s}.
\end{equation}
Finally, it remains only to control the  term  $\Sigma_3^5(s)$. 
Also, by using the fact that $\p\le C \q$, we have
\begin{eqnarray}\label{16nov4}
\Sigma_3^5(s)
 &\le& C\int_{B}\big(|\grad w|^2-(y.\grad  w)^2\big) \q{\mathrm{d}}y\\ 
&&\underbrace{-\b  \int_{B}|\grad  w|^2(1-|y|^2)\q\log (1-|y|^2) {\mathrm{d}}y}_{\chi_7(s)},
\quad  \forall  s\ge \max  (-\log T^*(x),1).\quad\quad\no
\end{eqnarray}
In a similar way to  the treatment of $\chi_1^2(s)$ and $\chi_2(s)$ to \eqref{X12} and \eqref{5nov11}, we write
\begin{equation}\label{6dec1}
\chi_7(s)\le C e^{-\frac{s}4} \int_{B}|\grad  w|^2 
\y+\frac{C}{s^{b-1}} \int_{B}|\grad  w|^2(1-|y|^2)\q. 
\end{equation}
By adding \eqref{16nov4} and \eqref{6dec1}, we get  for all  $s\ge \max  (-\log T^*(x),1)$
\begin{equation}\label{17nov15}
\Sigma_3^5(s) \le C\int_{B}\big(|\grad w|^2-(y.\grad  w)^2\big) \q
\y + C e^{-\frac{s}{4} }\int_{B}|\grad w|^2\y.
\end{equation}
Now, we are able to conclude the proof of the inequality
(\ref{cor17nov}). For this, we combine  \eqref{14nov5}, \eqref{14nova1}, \eqref{5dec6}, \eqref{17nov3}, \eqref{17nov4}   and \eqref{17nov15}  to get the desired estimate (\ref{cor17nov}) 
 which ends the proof of Lemma 2.4. \Box


\subsection{Existence of a decreasing functional for equation \eqref{var}}
In this subsection, by using 
Lemmas 2.3 and 2.4, 
 we are going to construct  a decreasing functional for  equation
$\eqref{var}$. Let us  define the following functional:
\begin{equation}\label{17novlyap1}
N(w(s),s)=exp \Big(\frac{p_c+3}{(b-1)s^{b-1}}\Big)K(w(s),s)+ \sigma e^{-\frac{\e_0}8s},
\end{equation}
where $\sigma$ is a constant will be determined later, $\e_0= 
 \frac{p_c-1}{32b(p_c+1)}$
and where
\begin{equation}\label{17novlyap}
K(w(s),s)=H(w(s),s)+ \frac{b}{2(s+s^{b+1})}L(w(s),s).
\end{equation}
We now state  the following proposition:
\begin{pro}\label{9avrilp}
 For all  $b\in (1,a)$ and $0<\e<   \frac{p_c-1}{32b(p_c+1)}$,
there exists $S_{3}\geq S_2$ and $\lambda_1>0$, such that for all $s \geq \max  (-\log T^*(x),S_{3})$, 
we have the following inequality:
\begin{eqnarray}\label{17nov55} 
N(w(s+1),s+1)-N(w(s),s)\!\!\!&\le &\!\!\!\!-\frac{\lambda_1}{s^{b}}\ia\iint\!\!
\Big( \frac{(\partial_{s}w)^2}{1-|y|^2}
+w^2+ |w|^{p_c+1}\Big)\ta
\y\t\qquad\quad\no\\
&&\!\!\!-\frac{\lambda_1}{s^b} \ia\iint (|\grad w|^2-(y.\grad w)^2)\ta\y\t.
\end{eqnarray}
Moreover, there exists
 $S_{4}\geq S_3$ such that for all $s \geq \max  (-\log T^*(x),S_{4})$, 
we have 
 \begin{equation}
N(w(s),s)\geq 0.
\end{equation}
\end{pro}
{\it Proof of Proposition 2.5}: Let  $b\in (1,a)$ and $0<\e<   \frac{p_c-1}{32b(p_c+1)}$.
 Combining the Lemmas 2.3  and 2.4, the fact that  $0\le \p\le \q$ and   choose $S_{3}\geq S_{2}$ large enough there exist
 $\mu_1>0$ such that 
 for all $s\geq \max  (-\log T^*(x),S_{3})$,
we have
\begin{eqnarray}\label{6dec11}
\frac{d}{ds}K(w(s),s)&\le &
\frac{p_c+3}{2s^{b}}K(w(s),s)- \frac{\mu_1}{s^b}\iint (|\grad w|^2-(y.\grad w)^2)\q\y\no\\
&&- \frac{\mu_1}{s^b}\iint \Big(\frac{(\partial_{s}w)^2}{1-|y|^2} +
  |w|^{p_c+1} + w^2\Big)\q \y\no\\
 &&+ C e^{-\frac{\e}{4} s}\int_{B}\Big(\frac{(\partial_{s}w)^2}{\sqrt{1-|y|^2}}  
+|\grad w|^2 +w^2\Big) \y+ Ce^{-s}.
\end{eqnarray}
Since,  for all $s\geq 1$, we have $C^{-1}\le exp \Big(\frac{p_c+3}{(b-1)s^{b-1}}\Big) \le C$. Then there exist $\mu_2>0$ such that  for all $s\geq \max  (-\log T^*(x),S_{3})$, we have
\begin{eqnarray}\label{7dec1}
\frac{d}{ds}N(w(s),s)&\le &- \frac{\mu_2}{s^b}\iint (|\grad w|^2-(y.\grad w)^2)\q\y\\
&&- \frac{\mu_2}{s^b}\iint \Big(\frac{(\partial_{s}w)^2}{1-|y|^2} +
  |w|^{p_c+1} + w^2\Big)\q \y\no\\
 &&+ C e^{-\frac{\e}{4} s}\int_{B}\Big(\frac{(\partial_{s}w)^2}{\sqrt{1-|y|^2}}  
+|\grad w|^2 +w^2\Big) \y+ Ce^{-s}-\frac{\e_0 \sigma}{8} e^{-\frac{\e_0 }{8}s}.\no
\end{eqnarray}
We now exploit the exponential space-time estimates \eqref{co2bis} and \eqref{co1bis}  in the particular case where $\eta =\frac{\e_0}{4(p_c+3)}$,  to show that
\begin{equation}\label{8dec1}
 C e^{-\frac{\e_0}{4} s}\ia\int_{B}\Big(\frac{(\partial_{s}w)^2}{\sqrt{1-|y|^2}}  
+|\grad w|^2 +w^2\Big) \y\t\le   C e^{-\frac{\e_0}{8} s}.
\end{equation}
By integrating the inequality \eqref{7dec1} in time between $s$ and $s+1$
 and taking into account $\e\in (0, \e_0)$, \eqref{8dec1} and choose $\sigma $ large enough  to deduce \eqref{17nov55}.

\medskip

To end the proof of the last point of Proposition 2.5, we refer the reader to 
 \cite{MZimrn05} and
\cite{omar1}. Let us mention that our proof strongly relies on the fact that $p<1+\frac{4}{N-2}$ which is implied by the fact that $p=p_{c}\equiv 1+\frac{4}{N-1}$.
\Box

\subsection{Proof of Theorem 1.1' }
 We define the following time:
\begin{equation}\label{new19dec1}
 t_1(x_0)=\max(T(x_0)-e^{-S_4},0).
\end{equation}
According to the Proposition 2.1, we obtain the following corollary which summarizes the principle properties of $N(w(s),s)$ defined in   \eqref{17novlyap1}.
\begin{cor}\label{19dec3} {\bf (Estimate on $N(w(s),s)$).}
For all $b\in (1,a)$,
  there exists  $t_1(x_0)\in [t_0(x_0), T(x_0))$
such that, for all $T_0\in (t_1(x_0),T(x_0)]$
 for all $s\ge  -\log (T_0-t_1(x_0))$
 and $x\in \er^N$ where $|x-x_0|\le \frac{e^{-s}}{\delta_0(x_0)}$,
 we have
$$-C\leq N(w(s),s) \leq N(w(\tilde{s}_{0}),\tilde{s}_{0}) ,$$
$$\int_{s}^{s+1}\iint\frac{(\partial_{s} w(y,\tau))^2}{1-|y|^2}\q{\mathrm{d}}y{\mathrm{d}}\tau\leq C\Big(1 +N(w(\tilde{s}_{0}),\tilde{s}_{0} ) \Big)s^{b},$$
$$ \int_{s}^{s+1}\int_{B_{1/2}}\Big(|\grad w(y,\tau)|^2+|w(y,\tau)|^{p_{c}+1}\Big)   {\mathrm{d}}y{\mathrm{d}}\tau\leq C\Big(1 +N(w(\tilde{s}_{0}),\tilde{s}_{0} ) \Big)s^{b},$$
where $\tilde{s}_{0}=-\log (T^*(x)-t_1(x_0))$.
\end{cor}

\medskip

\noindent  {\bf{Remark 2.2.}} Using the definition of  (\ref{scaling}) of
$w_{x,T^*(x)}=w$, we write easily
\begin{eqnarray}\label{cor0003}
   N(w(\widetilde{s_0}), \widetilde{s_0})
\le \widetilde{K_0},\nonumber
\end{eqnarray}
where
$\widetilde{K_0}=\widetilde{K_0}(T(x_0)-t_1(x_0),\|(u(t_1(x_0)),\partial_tu(t_1(x_0)))\|_{H^{1}\times L^{2}(B(x_0,\frac{T(x_0)-t_1(x_0)}{\delta_0(x_0)}))})$.

\medskip

 With    Corollary \ref{19dec3}, we are in a position to state and prove Theorem 1.1', which is a uniform
 version of Theorem 1.1 for $x$ near $x_0$.

\medskip

\noindent {\bf{Theorem} 1.1'{\it
{{\bf (Uniform polynomially space-time  estimate of  solution of  \eqref{A})}\label{t1bis}.\\
Let  $u $    a solution of ({\ref{gen}})
with blow-up graph $\Gamma:\{x\mapsto T(x)\}$ and  $x_0$ is a non
characteristic point. Then for all $b\in (1,a)$,  
for all $T_0\in [t_1(x_0), T(x_0)]$, for all
 $s\ge -\log (T_0-t_1(x_0))$
  and $x\in \er^N$, where $|x-x_0|\le \frac{e^{-s}}{\delta_0(x_0)}$,
we have
\begin{equation}\label{a1}
\int_{s}^{s+1} \!\!\int_{B}\!\!\Big((\partial_s
w(y,\tau))^2+|\grad
w(y,\tau)|^2
+|w (y,\tau)|^{p_c+1}\Big){\mathrm{d}}y{\mathrm{d}}\tau\le
K_2s^b,
\end{equation}
where $w=w_{x,T^*(x)}$,  where $T^*(x)$ is defined \eqref{18dec1} and
where the constant $K_2$ depends only on $N,  p, M, K_1, \delta_{0}(x_{0})$, $ T(x_{0})$ and 
$\|(u(t_1(x_0)),\partial_tu(t_1(x_0)))\|_{
H^{1}\times
L^{2}(B(x_0,\frac{T(x_0)-t_1(x_0)} {\delta_0(x_0)}) )}$.}}}

\bigskip

{{\it {Proof of Theorem 1.1'}}}:
Note that the estimate on the space-time $L^{2}$ norm of $\partial_{s} w$ was already proved in Corollary 2.5. Thus we focus on the space-time $L^{p_{c}+1}$ norm of $w$ and  $L^{2}$ norm of $\grad w$. This estimate was proved in Corollary 2.6 but just for the space-time $L^{p_{c}+1}$ norm of $w$ and  $L^{2}$ norm of $\grad w$ in $B_{1/2}$. To extend this estimate from $B_{1/2}$ to $B$ we refer the reader to Merle and Zaag \cite{MZimrn05} (unperturbed case) and Hamza and Zaag 
\cite{HZjhde12} 
 (perturbed case), where they introduce a new covering argument to extend the estimate of any known space $L^{q}$ norm of $w$, $\partial_{s} w$ or $\grad w$, from $B_{1/2}$ to $B$.
\Box
\section{Proof of Theorem \ref{t2}}
This section is devoted to the proof of
Theorem \ref{t2}. 
Throughout this section, we  assume $a>2$.
This section is divided into two parts:
\begin{itemize}
\item  In subsection 3.1,   based upon Theorem 1.1', we construct  a Lyapunov functional 
for equation (\ref{A}) and a blow-up criterion involving this functional.
\item In subsection 3.2, we prove   Theorem \ref{t2}.
\end{itemize}
\subsection{Existence  of a Lyapunov functional for equation (\ref{A}) and a blow-up criterion}
 Consider $u $   a solution of (\ref{gen}) with
blow-up graph $\Gamma:\{x\mapsto T(x)\}$ and  $x_0$ is a non
characteristic point.
Let   $T_0\in (t_1(x_0), T(x_0)]$, for all 
$x\in \er^N$  such that $|x-x_0|\le \frac{T_0}{\delta_0(x_0)}$, then we write $w$ instead of $w_{x,T^*(x)}$ defined in (\ref{scaling}) with $T^*(x)$ given in  (\ref{18dec1}).
Firstly,  for all
 $s\ge -\log (T_0-t_1(x_0))$  we introduce the following natural 
functional:
\begin{eqnarray}
E_0(w(s),s)\!\!\!&=&\!\!\!\!\iint \Big(\frac{1}{2}(\partial_{s}w)^2+\frac{1}{2}(|\grad w|^2-(y.\grad w)^2)+\frac{p_c+1}{(p_c-1)^2}w^2-\frac{|w|^{p_c+1}}{p_c+1}\Big) \y\quad\nonumber\\
&&-e^{\frac{-2(p_c+1)s}{p_c-1}} \iint  F(e^{\frac{2s}{p_c-1}}w) \y.\label{10dec1}
\end{eqnarray}
Moreover,  for all
 $s\ge -\log (T_0-t_1(x_0))$,  we define the functional
\begin{equation}\label{10dec2}
H_{0}(w(s),s)=E_{0}(w(s),s)+\frac{1 }{s^{\frac{a-2}{4}}}.
\end{equation}
We derive that the functional  $H_{0}(w(s),s)$ is a decreasing 
  functional  of time  for equation (\ref{A}),  provided that $s$ large enough.
Let us first control  the time derivative of the  functional $E_0(w(s),s)$ in the following lemma:
\begin{lem}\label{l10dec}
\label{energylyap0} For all $s\ge \max(-\log T^*(x_0),1)$, we have  the following inequality
\begin{equation}\label{10dec2}
 \frac{d}{ds}E_0(w(s),s)= - \int_{\partial B} (\partial_s w)^2
{\mathrm{d}}\sigma+\Sigma_{4}(s),
\end{equation}
where $\Sigma_{4}(s)$ satisfies
\begin{equation}\label{10dec3}
\Sigma_{4}(s)\leq \frac{C}{s^a}\iint |w|^{p_c+1}\y+Ce^{-s}.
\end{equation}
\end{lem} 
{\it Proof}: Multiplying \eqref{A} by $\partial_{s}w$ and  integrating over $B$, we obtain $\eqref{10dec2}$ where
\begin{equation}\label{nn1}
\Sigma_{4}(s)=\frac{2(p_c+1)}{p_c-1}e^{-\frac{2(p_c+1)s}{p_c-1}}\iint  F(e^{\frac{2s}{p_c-1}}w) \y-\frac{2}{p_c-1}e^{-\frac{2p_cs}{p_c-1}}\iint wf(e^{\frac{2s}{p_c-1}}w)\y.\\
\end{equation}
According to $\eqref{F1nov}$ and $\eqref{nn1}$, we get the desired estimate in \eqref{10dec3}, which end the proof of Lemma \ref{l10dec}.
\Box

\medskip 

With Lemma 3.1 and Theorem 1.1' we are in position to prove that $H_{0}(w(s),s)$ is a Lyapunov functional of equation $\eqref{A}$,
provided that $s$ is large enough.
 \begin{pro}
 Consider $u $   a solution of ({\ref{gen}}) with
blow-up graph $\Gamma:\{x\mapsto T(x)\}$ and  $x_0$ is a non
characteristic point.
Then, there exists $t_2(x_0)\in [t_1(x_0),T(x_0))$ such that,  for all   $T_0\in (t_2(x_0), T(x_0)]$, for all $s\ge -\log (T_0-t_2(x_0))$ and
$x\in \er^N$  such that $|x-x_0|\le \frac{e^{-s}}{\delta_0(x_0)}$,  we have
\begin{equation}\label{lyap} 
 H_{0}(w(s+1),s+1)- H_{0}(w(s),s) \leq -\int_{s}^{s+1}\int_{\partial B} (\partial_{s}w(\sigma ,\tau))^2{\mathrm{d}}\sigma \t.
\end{equation}
Moreover, there exists $S_{6}\geq S_{4}$ such that, for all $s\geq max(
-\log (T_0-t_2(x_0)),S_6)$, we have: $$H_{0}(w(s),s)\geq 0.$$
\end{pro}
{\it Proof}: We apply  the polynomial space-time estimate \eqref{a1} in the particular case where $b=\frac{a}2>1$
 and Lemma 3.1 to get,  for all  $s\ge -\log (T_0-t_1(x_0))$
\begin{equation}\label{11dec1}
 E_{0}(w(s+1),s+1)-E_{0}(w(s),s) \leq 
-\int_{s}^{s+1}\int_{\partial B} (\partial_{s}w(\sigma ,\tau))^2{\mathrm{d}}\sigma \t  +\frac{C}{s^{\frac{a}2}}.
\end{equation}
Then we write,    for all  $s\ge -\log (T_0-t_1(x_0))$
\begin{eqnarray}\label{11dec2}
 H_{0}(w(s+1),s+1)-H_{0}(w(s),s) &\leq &
-\int_{s}^{s+1}\int_{\partial B} (\partial_{s}w(\sigma ,\tau))^2{\mathrm{d}}\sigma \t\no\\
&&  +\frac{C}{s^{\frac{a}2}}+\frac{1}{(s+1)^{\frac{a-2}4}}-\frac{1}{s^{\frac{a-2}4}}.
\end{eqnarray}
For all $s\ge 1$, by the mean value theorem to the function $x\longmapsto \frac{1}{x^{\frac{a+2}{4}}}$, 
 between $s$ and $s+1$, so we can say that there exists a constant $\gamma=\gamma(s)\in ( 0,1) $ such that
\begin{equation}\label{11dec1}
 \frac{1 }{(s+1)^{\frac{a-2}{4}}}-\frac{1 }{s^{\frac{a-2}{4}}}=
\frac{2-a}{4(s+\gamma )^{\frac{a+2}{4}}}<\frac{2-a}{4(s+1)^{\frac{a+2}4}}\le\frac{2-a}{2^{\frac{a+10}4}}\frac1{s^{\frac{a+2}{4}}}.
 \end{equation}
Finally, by exploiting \eqref{11dec2}  and the inequality \eqref{11dec1},
we can 
 choose $S_5>S_4$  large enough so that we  get \eqref{lyap}, where 
\begin{equation}\label{3jan1}
 t_2(x_0)=\max(T(x_0)-e^{-S_5},0).
\end{equation}
To end the proof of the last point of Proposition 3.2, we refer the reader to  \cite{MZimrn05}. 
This concludes the proof of Proposition 3.2.
\Box
\subsection{Proof of Theorem \ref{t2} }
In this subsection, we prove Theorem \ref{t2}. Note that the lower bound follows from the finite speed of propagation and the wellposedness in $H^1\times L^2$. For a detailed argument in the similar case of equation \eqref{A}, see Lemma 3.1 p 1136 in  \cite{MZimrn05}.
Let us first use Proposition 4 
and the averaging technique of \cite{MZimrn05} and \cite{MZma05} to get the following bounds:
\begin{cor}\label{fin}
{\rm For all} $s \geq \max( -\log (T_0-t_2(x_0))
,S_6)$,
{\rm it holds that}
\begin{eqnarray*}
-C \leq E_{0}(w(s),s)&\leq &M_0,\\
 \int_{s}^{s+1}\int_{\partial B}(\partial_{s}w(\sigma,\tau))^2\t &\leq &M_0,\\
\int_{s}^{s+1}\int_{ B}  \Big(\partial_{s}w(y,\tau )-\lambda (\tau ,s)w(y,\tau )\Big)^2dyd\tau &\leq &CM_0,
\end{eqnarray*}
where $0\leq\lambda (\tau ,s)\leq C$.
\end{cor}

Proof of Theorem \ref{t2}: The proof  is similar to the one in the unperturbed case treated by Merle and Zaag in \cite{MZimrn05} and \cite{MZma05}   and also used by Hamza and Zaag in \cite{HZjhde12}, 
\cite{HZnonl12} and Hamza and Saidi in \cite{omar1} and \cite{omar2}. To be accurate and concise in our results, there is an analogy between the exponential smallness exploited in \cite{HZjhde12} by Hamza and Zaag and the polynomial smallness used here. The unique difference lies in the treatment of the perturbed term which is treated by Hamza and Saidi \cite{omar1} and \cite{omar2}.  This concludes the proof of Theorem \ref{t2}.
\Box



\appendix
\section{Some identity related to the Pohozaev multiplier }
In this appendix,  for all $s\ge \max (-\log T^*(x),1)$, we evaluate  the term
\begin{equation}\label{A0}
{\cal{L}}(s)=\int_{B}(y.\grad  w) \Big(\div(\nabla w- (y.\nabla w)y)\Big)\p \y,
\end{equation}
where \begin{equation}\label{A50}
\p=(1-|y|^2)^{1+\frac1{s^b}}\big(1-\log (1-|y|^2)\big).
\end{equation}
  More precisely, we prove 
the following identity:
\begin{lem}\label{A540} 
For all $w\in {\cal  H}$ it holds that
\begin{eqnarray}\label{import1}
{\cal{L}}(s)  &=&
(1+\b)\int_{B}|\grad_{\theta} w|^2|y|^2\q\log(1-|y|^2)
{\mathrm{d}}y-
\b
\int_{B}|\grad_{\theta} w|^2|y|^2\q
{\mathrm{d}}y\no\\
&&+\frac{N-2}2\int_{B}\big(|\grad w|^2-(y.\grad  w)^2\big) \p {\mathrm{d}}y \\
&&+\b  \int_{B}(y.\grad  w)^2\p {\mathrm{d}}y- \int_{B}(1-|y|^2)(y.\grad  w)^2
\q {\mathrm{d}}y.\no
\end{eqnarray}
\end{lem}
{\it Proof}: 
We divide ${\cal{L}}(s)$ into two terms:
 $ {\cal{L}}(s)={\cal{L}}_1(s)+{\cal{L}}_2(s),$ where
\begin{equation}\label{A1}
{\cal{L}}_1(s)=\int_{B}(y.\grad  w) \Delta w\p \y,
\end{equation}
and
\begin{equation}\label{A2}
{\cal{L}}_2(s)=-\int_{B}(y.\grad  w) \div( (y.\nabla w)y)\p \y.
\end{equation}
To estimate ${\cal{L}}_1(s)$, we start
observe  the immediate  identity
\begin{equation}\label{A3}
  (y.\grad w) \Delta w= \sum_{i,j}y_i\partial_i w
\partial^2_jw.
\end{equation}
By  integrating by parts, exploiting \eqref{A3} and  the fact that  $\ds{\sum_{i,j}\delta_{i,j}\partial_i w
\partial_jw=|\grad w|^2}$,      we can write
\begin{eqnarray}\label{A4}
{\cal{L}}_1(s) &=&-\frac12\sum_{i,j}\int_{B}y_i\partial_i((
\partial_j w)^2)\p \y-\int_{B}|\grad w|^2
\p\y \no\\ &&-\sum_{i,j}\int_{B}y_i\partial_i w
\partial_j w\partial_j\p \y.
\end{eqnarray}
By using the identity  $$\partial_j\p=-2(1+\b)\frac{y_j}{1-|y|^2}\p+2y_j\q,$$ 
\eqref{A4} and  integrating by part one has that
\begin{eqnarray}\label{A5}
{\cal{L}}_1(s)  &=&\frac12\int_{B}|\grad w|^2\div (\p y)
{\mathrm{d}}y-\int_{B}|\grad w|^2 \p
{\mathrm{d}}y\no\\
&&+2(1+\b)\int_{B}(y.\grad w)^2\frac{\p}{1-|y|^2}
 \y-2\int_{B}(y.\grad w)^2\q\y.
\end{eqnarray}
Furthermore, by using the identity \eqref{12nov3} and  \eqref{A5}, we get 
\begin{eqnarray}\label{A6}
{\cal{L}}_1(s)  &=&-(1+\b)\int_{B}|\grad w|^2\frac{|y|^2\p}{1-|y|^2}
{\mathrm{d}}y+2(1+\b)\int_{B}(y.\grad w)^2\frac{\p}{1-|y|^2}\y\\
&&+\frac{N-2}2\int_{B}|\grad w|^2 \p
{\mathrm{d}}y 
 -2\int_{B}(y.\grad w)^2\q\y+\int_{B}|y|^2|\grad w|^2\q\y.\no
\end{eqnarray}
To estimate ${\cal{L}}_2(s)$, we start use the classical identity
\begin{equation}\label{A11}
\div( (y.\nabla w)y)= N (y.\nabla w)+ \grad (y.\nabla w).y,
\end{equation}
and  integrating by part, to obtain
\begin{equation}\label{A12}
{\cal{L}}_2(s) =-N\int_{B}(y.\grad
w)^2\p{\mathrm{d}}y +\frac12 \int_{B}(y.\grad w)^2 \div (
\p y) {\mathrm{d}}y.
\end{equation}
Also, by using \eqref{12nov3},
we can write
\begin{eqnarray}\label{A13}
{\cal{L}}_2(s)&=&-\frac{N}2\int_{B}(y.\grad  w)^2\p
{\mathrm{d}}y-(1+\b ) \int_{B}(y.\grad  w)^2
\frac{|y|^2\p}{1-|y|^2} {\mathrm{d}}y\no\\
&&+ \int_{B}|y|^2(y.\grad  w)^2
\q {\mathrm{d}}y.
\end{eqnarray}
By combining  \eqref{A6} and  \eqref{A13} and using the expression of $\p$ defined in \eqref{A50} to write the basic idendity   
$\ds{\frac{\p}{1-|y|^2} -
\q=-\q \log (1-|y|^2)}$,
we deduce easily
\eqref{import1}, which ends the proof of Lemma \ref{A540}. 
\Box
\def\cprime{$'$} \def\cprime{$'$}
\providecommand{\bysame}{\leavevmode\hbox to3em{\hrulefill}\thinspace}
\providecommand{\MR}{\relax\ifhmode\unskip\space\fi MR }
\providecommand{\MRhref}[2]{%
  \href{http://www.ams.org/mathscinet-getitem?mr=#1}{#2}
}
\providecommand{\href}[2]{#2}


\noindent{\bf Address}:\\
Universit\'e de Tunis El Manar, Facult\'e des Sciences de Tunis, LR03ES04 \'Equations aux d\'eriv\'ees partielles et applications, 2092 Tunis, Tunisie\\
\vspace{-7mm}
\begin{verbatim}
e-mail: ma.hamza@fst.rnu.tn
\end{verbatim}
\end{document}